\documentclass{svmult}

\usepackage{graphicx}    
\usepackage{multicol}    

\usepackage{amsmath,amssymb}
\def\R{{\mathbb{R}}}

\newcommand{\bvec}{\mathbf{B}}

\newcommand{\evec}{\mathbf{E}}

\newcommand{\fracp}[2]{\frac{\partial #1}{\partial #2}} 
\renewcommand{\div}{\nabla\cdot}
\newcommand{\curl}{\nabla\times}
\newcommand{\eps}{\epsilon}           

\newcommand{\beq}{\begin{equation}}
\newcommand{\eeq}{\end{equation}}

\begin{document}
\title*{An adaptive numerical method for the Vlasov equation
based on a multiresolution analysis}
\titlerunning{An adaptive numerical method for the Vlasov equation}
\author{
N. Besse\inst{1}, 
F. Filbet\inst{2}, 
M. Gutnic\inst{2}, I. Paun \inst{2},
E. Sonnendr{\"u}cker\inst{2}}
\institute{C.E.A, BP 12, 91680 Bruy{\`e}res-le-Ch{\^a}tel, France,
\textit{nicolas.besse@cea.fr}
\and IRMA, Universit{\'e} Louis Pasteur, 67084 Strasbourg cedex, France,
\textit{filbet,gutnic,ipaun,sonnen@math.u-strasbg.fr}}

\maketitle
\bigskip
\section{Introduction}
\label{sec:intro}
Plasmas, which are gases of charged particles, and
charged particle beams can be described by a 
distribution function $f(t,x,v)$ dependent on time
$t$, on position $x$ and on velocity $v$. The function $f$
represents the probability of presence of a particle at position
$(x,v)$ in phase space at time $t$.
It satisfies the so-called Vlasov equation
\begin{equation}
\label{vlasov} 
\frac{\partial f}{\partial t} + v\cdot\nabla_x f + F(t,x,v)\cdot\nabla_v f = 0.
\end{equation}
The force field $F(t,x,v)$ consists of applied and self-consistent
electric and magnetic fields: 
$$
F=\frac{q}{m}(E_{self}+E_{app}+v\times(B_{self}+B_{app})), $$
where $m$ represents the mass of a particle and $q$ its charge.
The self-consistent part of the force field is solution of Maxwell's equations
\begin{eqnarray*}
-\frac{1}{c^2} \fracp{\evec}{t}+\curl\bvec = \mu_0{\mathbf j}, &&\hspace{1cm}
\div\evec = \frac{\rho}{\eps_0}, \\
 \fracp{\bvec}{t}+\curl\evec = 0,&&\hspace{1cm}
\div\bvec = 0.
\end{eqnarray*}  

The coupling with the Vlasov equation results from the source terms
$\rho$ and ${\mathbf j}$ such that: 
$$
\rho(t,x) = q \int_{\R^d} f(t,x,v)\,dv,  ~~~ 
{\mathbf j}=q \int_{\R^d} f(t,x,v)v \,dv  .
$$
We then obtain the nonlinear Vlasov-Maxwell equations.
In some cases, when the field are slowly varying the magnetic field becomes
negligible and the Maxwell equations can be replaced by the Poisson equation
where:
\begin{equation}
\label{poisson}
E_{self}(t,x) = -\nabla_x \phi(t,x), \quad -\varepsilon_0 \Delta_x\phi = \rho.
\end{equation}

The numerical resolution of the Vlasov equation is usually performed
by particle methods (PIC) which consist in approximating the plasma by
a finite number of particles. The trajectories of these particles are
computed from the characteristic curves given by the Vlasov equation,
whereas self-consistent fields are computed on a mesh of the physical
space. This method allows to obtain satisfying results with a few
number of particles. However, it is well known that, in some cases,
the numerical noise inherent to the particle method becomes too
important to have an accurate description of the distribution function
in phase space. Moreover, the numerical noise only decreases in $\sqrt
N$, when the number of particles $N$ is increased.  To remedy to this
problem, methods discretizing the Vlasov equation on a mesh of phase
space have been proposed. A review of the main methods for the
resolution of the Vlasov equation is given in these proceedings
\cite{sonnen:enumath}.

The major drawback of methods using a uniform and fixed mesh is that
their numerical cost is high, which makes them rather inefficient when
the dimension of phase-space grows.  For this reason we are
investigating here a method using an adaptive mesh.  The adaptive
method is overlayed to a classical semi-Lagrangian method which is
based on the conservation of the distribution function along
characteristics. Indeed, this method uses two steps to update the value
of the distribution function at a given mesh point. The first one consists in
following the characteristic ending at this mesh point backward in
time, and the second one in interpolating its value there
from the old values at the surrounding mesh points. Using the conservation
of the distribution function along the characteristics this will yield
its new value at the given mesh point.
This idea was originally introduced by Cheng and Knorr
\cite{cheng} along with a time splitting technique enabling to compute
exactly the origin of the characteristics at each fractional step. In
the original method, the interpolation was performed using cubic
splines. This method has since been used extensively by plasma
physicists (see for example \cite{feix,ghizzo} and the references
therein).  It has then been generalized to the frame of semi-Lagrangian
methods by E.  Sonnendr{\"u}cker {\em et al.}
\cite{sonnen}. This method has also been used to investigate problems
linked to the propagation of strongly nonlinear  heavy ion beams \cite{hif}.

In the present work, we have chosen to introduce a phase-space mesh
which can be refined or derefined adaptively in time.
For this purpose, we use a technique based on multiresolution analysis
which is in the same spirit as the methods developed in particular  by
S. Bertoluzza \cite{bertoluzza}, A. Cohen {\em et al.} \cite{Cal2001} 
and M. Griebel and F. Koster
\cite{GK2000}. 
We represent the distribution function on a wavelet basis at different
scales.  We can then compress it by eliminating coefficients which are
small and accordingly remove the associated mesh points. Another
specific feature of our method is that we use an advection in physical
and velocity space forward in time to predict the useful grid points
for the next time step, rather than restrict ourselves to the
neighboring points. This enables us to use a much larger time step, as
in the semi-Lagrangian method the time step is not limited by a
Courant condition. Once the new mesh is predicted, the semi-Lagrangian
methodology is used to compute the new values of the distribution
function at the predicted mesh points, using an interpolation based
on the wavelet decomposition of the old distribution function.
The mesh is then refined again by performing a wavelet transform,
and eliminating the points associated to small coefficients.

This paper is organized as follows.
In section \ref{sec:MRA}, we recall the tools of multiresolution analysis
which will be needed for our method, precizing what kind of wavelets
seem to be the most appropriate in our case.
Then, we describe in section
\ref{sec:algo} the algorithm used in our method, first for 
the non adaptive mesh case and then for the adaptive mesh case.
Finally we present a few preliminary numerical results.

\section{Multiresolution analysis}
\label{sec:MRA}

The semi-Lagrangian method consists mainly of two steps, an advection
step and an interpolation step. The interpolation part is performed
using for example a Lagrange interpolating polynomial on a uniform grid.
Thus interpolating wavelets provide a natural way to extend this procedure
to an adaptive grid in the way we shall now shortly describe.
 
For simplicity, we shall restrict our description to the 1D case of
the whole real line. It is straightforward to extend it to periodic boundary
conditions and it can also be extended to an interval with Dirichlet
boundary conditions. The extension to higher dimension is performed
using a tensor product of wavelets and will be addressed at the end of
the section.
 
For any value of $j\in\mathbb{Z}$, we consider a uniform grid $G^j$
of step $2^{-j}$. The grid points are located at $x^j_k=k2^{-j}$.
This defines an infinite sequence of grids that we denote by
$(G_j)_{j\in\mathbb{Z}}$, and $j$ will be called the level of the grid.

In order to go from one level to the next or the previous, we define
a projection operator and a prediction operator. Consider two grid
levels $G_j$ and $G_{j+1}$ and discrete values (of a function)
denoted by $(c_k^j)_{k\in\mathbb{Z}}$ and $(c_k^{j+1})_{k\in\mathbb{Z}}$.
Even though we use the same index $k$ for the grid points in the two
cases, there are of course twice as many points in any given interval
on $G_{j+1}$ as on $G_j$. 
Using the terminology in \cite{Cal2001}, we then define the projection operator
\begin{align*}
P_{j+1}^j : G_{j+1} &\rightarrow G_j, \\
c_{2k}^{j+1} &\mapsto c_k^j,
\end{align*}
which is merely a restriction operator, as well as the prediction operator
\begin{align*}
P_j^{j+1} : &\;G_j \rightarrow G_{j+1}, \\
\mbox{ such that } &c_{2k}^{j+1} = c_k^j, \\
& c_{2k+1}^{j+1} = P_{2N+1}(x_{2k+1}^{j+1}), \\
\end{align*}
where $P_{2N+1}$ stands for the Lagrange interpolation
polynomial of odd degree $2N+1$ centered at the point $(x_{2k+1}^{j+1})$.

Using the just defined prediction operator, we can construct on $G_j$
a subspace of $L^2(\mathbb{R})$ that we shall denote by $V_j$, a basis
of which being given by $(\varphi_k^j)_{k\in\mathbb{Z}}$ such that
$\varphi_k^j(x_{k'}^j)=\delta_{kk'}$ where $\delta_{kk'}$is the
Kronecker symbol. The value of $\varphi_k^j$ at any point of the real
line is then obtained by applying, possibly an infinite number of
times, the prediction operator.

In the wavelets terminology the $\varphi_k^j$ are called scaling functions. 
We shall also denote by
$\varphi=\varphi_0^0$. Let us notice that
$$\varphi_k^j(x)=\varphi(2^jx-k).$$
It can be easily verified that the scaling functions satisfy the following 
properties:
\begin{itemize}
\item Compact support: the support of $\varphi$ is included in
$[-2N-1,2N+1]$.
\item Interpolation: by construction $\varphi(x)$ is interpolating
in the sense that $\varphi(0)=1$ and $\varphi(k)=0$ if $k\neq 0$.
\item Polynomial representation: all polynomials of degree less or equal to
$2N+1$ can be expressed exactly as linear combinations of the $\varphi_k^j$.
\item Change of scale: the $\varphi$ at a given scale can be expressed
as a linear combination of the $\varphi$ at the scale immediately below:
$$\varphi(x)=\sum_{-2N-1}^{2N+1} h_l\varphi(2x-l).$$
\end{itemize}

Moreover the sequence of spaces $(V_j)_{j\in\mathbb{Z}}$ defines
a multiresolution analysis of
$L^2(\mathbb R)$, i.e. it satisfies the following properties:
\begin{itemize}
\item $\ldots\subset V_{-1}\subset V_{0}\subset V_{1}\subset \ldots
\subset  V_{n}\subset \ldots \subset L^2(\mathbb R)$. 
\item $\cap V_j=\{0\}$, $\overline{\cup V_j}=L^2(\mathbb R)$.
\item $f\in V_j \leftrightarrow f(2~\cdot) V_{j+1}$.
\item $\exists\, \varphi$ (scaling function) such that
$\{\varphi(x-k)\}_{k\in{\mathbb Z}}$ 
is a basis of $V_0$ and
$\{\varphi^j_k=2^{j/2}\varphi(2^j\,x-k)\}_{k\in{\mathbb Z}}$ 
is a basis of $V_j$.
\end{itemize}

As $V_j\subset V_{j+1}$, there exists a supplementary of
$V_j$ in $V_{j+1}$ that we shall call the detail space and
denote by $W_j$ :
$$
V_{j+1}=V_{j}\oplus W_j.
$$
The construction of $W_j$ can be made in the following way: an element
of $V_{j+1}$ is characterized by the sequence$(c_k^{j+1})_{k\in{\mathbb Z}}$
and by construction we have $c_k^j=c_{2k}^{j+1}$. Thus, if we define
$d_k^j=c_{2k+1}^{j+1}-P_{2N+1}(x_{2k+1}^{j+1})$, where $P_{2N+1}$ is the
Lagrange interpolation polynomial by which the value of an element of $V_j$ at 
the point $(x_{2k+1}^{j+1})$ can be computed, $d_k^j$
represents exactly the difference between the value in $V_{j+1}$
and the value predicted in $V_j$. Finally, any
element of $V_{j+1}$ can be characterized by the two sequences
$(c_k^j)_k$ of values in $V_j$ and $(d_k^j)_k$ of details in
$W_j$. Moreover this strategy for constructing $W_j$ is particularly
interesting for adaptive refinement as $d_k^j$ will be small at places where
the prediction from $V_j$ is good and large elsewhere, which gives us a natural
refinement criterion.
Besides, there exists a function $\psi$, called wavelet such that
$\{\psi^j_k=2^{j/2}\psi(2^j\,x-k)\}_{k\in{\mathbb Z}}$
is a basis of $W_j$.

In practise, for adaptive refinement we set the coarsest level $j_0$
and the finest level $j_1$, $j_0<j_1$, and we decompose the space
corresponding to the finest level on all the levels in between: 
$$
V_{j_1}=V_{j_0}\oplus W_{j_0}\oplus W_{j_0+1}\oplus \cdots \oplus
W_{j_1-1}.  $$ 
A function $f\in V_{j_1}$ can then be decomposed as follows 
$$
f(x)=\sum_{l=-\infty}^{+\infty}c_l^{j_0}\,\varphi_l^{j_0}(x)
+\sum_{j=j_0}^{j_1-1}\sum_{l=-\infty}^{+\infty} d_l^j\,\psi_l^j(x), $$
where the $(c_l^{j_0})_l$ are the coefficients on the coarse mesh and
the $(d_l^j)_l$ the details at the different level in between.
\setlength{\unitlength}{4144sp}%
\begingroup\makeatletter\ifx\SetFigFont\undefined%
\gdef\SetFigFont#1#2#3#4#5{%
  \reset@font\fontsize{#1}{#2pt}%
  \fontfamily{#3}\fontseries{#4}\fontshape{#5}%
  \selectfont}%
\fi\endgroup%
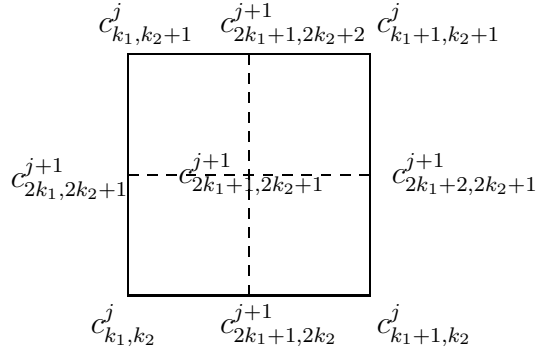
\begin{figure}
\centerline{
\begin{picture}(1890,1980)(2836,-1861)
\thinlines
{\put(3151,-1636){\framebox(1440,1440){}}
}%
{\multiput(3871,-196)(0.00000,-115.20000){13}{\line( 0,-1){ 57.600}}
}%
{\multiput(3151,-916)(115.20000,0.00000){13}{\line( 1, 0){ 57.600}}
}%
\put(4726,-961){\makebox(0,0)[lb]{\smash{\SetFigFont{12}{14.4}{\rmdefault}{\mddefault}{\updefault}{$c_{2k_1+2,2k_2+1}^{j+1}$}%
}}}
\put(2926,-1861){\makebox(0,0)[lb]{\smash{\SetFigFont{12}{14.4}{\rmdefault}{\mddefault}{\updefault}{$c_{k_1,k_2}^j$}%
}}}
\put(2448,-991){\makebox(0,0)[lb]{\smash{\SetFigFont{12}{14.4}{\rmdefault}{\mddefault}{\updefault}{$c_{2k_1,2k_2+1}^{j+1}$}%
}}}
\put(3450,-961){\makebox(0,0)[lb]{\smash{\SetFigFont{12}{14.4}{\rmdefault}{\mddefault}{\updefault}{$c_{2k_1+1,2k_2+1}^{j+1}$}%
}}}
\put(3700,-61){\makebox(0,0)[lb]{\smash{\SetFigFont{12}{14.4}{\rmdefault}{\mddefault}{\updefault}{$c_{2k_1+1,2k_2+2}^{j+1}$}%
}}}
\put(2971,-61){\makebox(0,0)[lb]{\smash{\SetFigFont{12}{14.4}{\rmdefault}{\mddefault}{\updefault}{$c_{k_1,k_2+1}^j$}%
}}}
\put(4636,-61){\makebox(0,0)[lb]{\smash{\SetFigFont{12}{14.4}{\rmdefault}{\mddefault}{\updefault}{$c_{k_1+1,k_2+1}^j$}%
}}}
\put(3700,-1861){\makebox(0,0)[lb]{\smash{\SetFigFont{12}{14.4}{\rmdefault}{\mddefault}{\updefault}{$c_{2k_1+1,2k_2}^{j+1}$}%
}}}
\put(4636,-1861){\makebox(0,0)[lb]{\smash{\SetFigFont{12}{14.4}{\rmdefault}{\mddefault}{\updefault}{$c_{k_1+1,k_2}^j$}%
}}}
\end{picture}}
\caption{\label{pred2d}Mesh refinement in  2D.}
\end{figure}

In two dimensions, the prediction operator which defines the
multiresolution analysis  is constructed by tensor product from the 1D operator.
In practise three different cases must be considered
 (see figure \ref{pred2d} for notations):
\begin{enumerate}
\item Refinement in $x$ (corresponding to points $c^{j+1}_{2k_1+1,2k_2}$
and $c^{j+1}_{2k_1+1,2k_2+2}$): we use the 1D prediction operator in
$x$ for fixed $k_2$.
\item Refinement in $v$ (corresponding to points $c^{j+1}_{2k_1,2k_2+1}$
and $c^{j+1}_{2k_1+2,2k_2+1}$): we use the 1D prediction operator in $v$ for 
fixed $k_1$.
\item Refinement in $v$ (corresponding to point $c^{j+1}_{2k_1+1,2k_2+1}$): 
we first  use the 1D prediction operator in $v$ for fixed  $k_1$
to determine the points which are necessary for applying the 1D prediction 
operator in
$x$ for fixed $k_2$ which we then apply.
\end{enumerate}
The corresponding wavelet bases are respectively of type
$\psi(x)\varphi(v)$, $\varphi(x)\psi(v)$ and $\psi(x)\psi(v)$ where
$\varphi$ and $\psi$ are respectively the scaling function and the
1D wavelet. We then obtain a 2D wavelet decomposition of the following form:
\begin{multline}\label{eq:wvltdec}
 f(x,v)=\sum_{k_1,k_2} \bigg(  c^{j_0}_{k_1,k_2} \, \varphi^{j_0}_{k_1}(x) 
\, \varphi^{j_0}_{k_2}(v)
     +\sum_{j_0}^{j_1-1} \Big(
                   d^{row,j}_{k_1,k_2} \,  \psi^{j}_{k_1}(x) \, 
\varphi^{j}_{k_2}(v)\\
                 + d^{col,j}_{k_1,k_2} \, \varphi^{j}_{k_1}(x) \,
                 \psi^{j}_{k_2}(v) + d^{mid,j}_{k_1,k_2} \,
                 \psi^{j}_{k_1}(x) \, \psi^{j}_{k_2}(v)\Big)\bigg).
\end{multline}

\section{The algorithms}
\label{sec:algo}

We want to numerically solve the Vlasov equation
(\ref{vlasov}) given an initial value of the
distribution function $f_0$. 

We start by describing the method based on an interpolation using the
wavelet decomposition of $f$ in the non adaptive case. Then we overlay an
adaptive algorithm to this method.

For those two algorithms, we first pick the resolution levels for the
phase-space meshes, from the coarsest $j_0$ to the finest
$j_1$. Although these levels could be different in $x$ and $v$, we
consider here for the sake of conciseness and clarity that they are
identical.

We also compute our scaling function on a very fine grid so that we can
obtain with enough precision its value at any point.

\subsection{The non adaptive algorithm}
\label{sec:NAalgo}

We are working in this case on the finest level corresponding to $j_1$
keeping all the points.

\bigskip

{\bf Initialization:} We decompose the initial condition in
the wavelet basis by computing the coefficients $c_{k_1,k_2}$ of the
decomposition in $V_{j_0}$ for the coarse mesh, and then adding the
details $d^j_{k_1,k_2}$ in the detail spaces $W_j$ for all the other
levels $j=j_0,\ldots,j_1-1$.  We then compute the initial electric
field.

\bigskip
{\bf Time iterations:}
\begin{itemize}
\item
  {\bf Advection in $x$:} We start by computing for each mesh point
  the origin of the corresponding characteristic exactly, the
  displacement being $v_j\Delta t$.  As we do not necessarily land on
  a mesh point, we compute the values of the distribution function at
  the intermediate time level, denoted by $f^*$, at the origin of the
  characteristics by interpolation from $f^n$. We use for this the
  wavelet decomposition (\ref{eq:wvltdec}) applied to $f^n$ from which
  we can compute $f^n$ at any point in phase space.

\item {\bf Computation of the electric field:}
We compute the charge density by integrating $f^*$ with respect to
$v$, then the electric field by solving the Poisson equation (this
step vanishes for the linear case of the rotating cylinder where the
advection field is exactly known).
\item{\bf Advection in $v$:} We start
by computing exactly the origin of the characteristic for each mesh
point, the displacement being $E(t^n,x_i)\Delta t$.  As we do not
necessarily land on a mesh point, we compute the values of the
distribution function at the intermediate time level, denoted by
$f^{n+1}$, at the origin of the characteristics by interpolation from
$f^*$. We use for this the wavelet decomposition of $f^*$ given by
(\ref{eq:wvltdec}) used at the previous step.
\end{itemize}

\subsection{The adaptive algorithm}
\label{sec:Aalgo}
In the initialization phase, we first compute the wavelet
decomposition of the initial condition $f_0$, and then proceed by
compressing it, i.e. eliminating the details which are
smaller than a threshold that we impose.  We then construct an
adaptive mesh which, from all the possible points at all the
levels between our coarsest and finest, contains only those of the coarsest
and those corresponding to details which are above the threshold. We
denote by $\tilde{G}$ this mesh.
\begin{itemize}
\item {\bf Prediction in $x$:}  We predict the 
positions of points where the details should be important at the next
time split step by advancing in $x$ the characteristics originating
from the points of the mesh $\tilde{G}$. For this we use an explicit
Euler scheme for the numerical integration of the
characteristics. Then we retain the grid points, at one level finer as
the starting point, surrounding the end point the characteristic.
\item {\bf Construction of mesh $\hat{G}$:} From the predicted mesh
$\tilde{G}$, we construct the mesh $\hat{G}$ where the values of the
distribution at the next time step shall be computed. This mesh
$\hat{G}$ contains exactly the points necessary for computing the
wavelet transform of $f^{*}$ at the points of $\tilde{G}$.
\item {\bf Advection in $x$:} As in the non adaptive case.
\item {\bf Wavelet transform of $f^{*}$:} We compute the
$c_k$ and $d_k$ coefficients at the  points of $\tilde{G}$ from the values of
$f^{*}$ at the points of  $\hat{G}$.
\item {\bf Compression:} We eliminate the points of $\tilde{G}$ where
the details $d_k$ are lower than the fixed threshold.
\item {\bf Computation of the electric field:} As in the non adaptive case.
\item {\bf Prediction in $v$:} As for the prediction in $x$.
\item {\bf Construction of mesh $\hat{G}$:} As previously.
This mesh $\hat{G}$ contains exactly the points necessary for computing the
wavelet transform of $f^{n+1}$ at the points of $\tilde{G}$ determined in
the prediction in $v$ step.
\item {\bf Advection in $v$:} As in the non adaptive case.
\item {\bf Wavelet transform of $f^{n+1}$:} We compute the
$c_k$ and $d_k$ at the points of $\tilde{G}$ from the values of
$f^{n+1}$ at the points of $\hat{G}$.
\item {\bf Compression:} We eliminate the points of $\tilde{G}$ where
the details $d_k$ are lower than the fixed threshold.
\end{itemize}
\section{Numerical results}
\label{sec:resnum}

We show here our first results obtained with the adaptive method.  We
consider first a linear problem, namely the test case of the rotating
cylinder introduced by Zalesak \cite{zalesak} to test advection
schemes. Then we consider a classical nonlinear Vlasov-Poisson test
case, namely the two stream instability.

\subsection{The slit rotating cylinder}
\label{sec:cyl}

We consider the following initial condition:
$$
f(0,x,v) = 
\left\{
\begin{array}{ll}
1 & \mbox{ if } \sqrt{x^2+v^2}<0.5 \mbox{ and if }x<0 \mbox{ or }
|v|>0.125,\\
0 & \mbox{ else}.
\end{array}\right.
$$ 
The computational domain is $[-0.5,0.5]\times [-0.5,0.5]$.

The advection field is $(v,-x)$, which corresponds to the Vlasov equation 
with an applied electric field $E_{app}(x,t)=-x$ and without self-consistent
field.
\begin{figure}
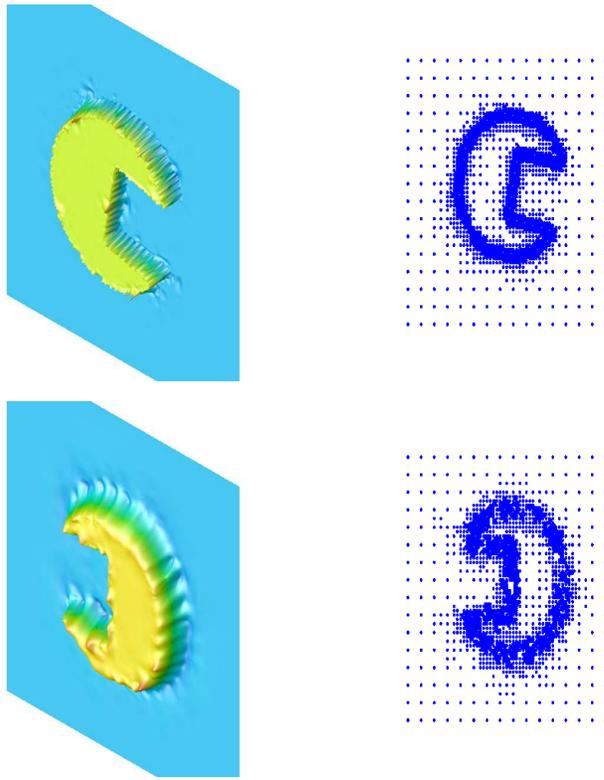

\begin{tabular}{cc}
\includegraphics[width=5cm,height=5.cm,angle=-90]{fxv.0.epsf} & 
\includegraphics[width=5cm,height=5.cm,angle=-90]{gridt.0.epsf}\\
\includegraphics[width=5cm,height=5.cm,angle=-90]{fxv.15.epsf} & 
\includegraphics[width=5cm,height=5.cm,angle=-90]{gridt.15.epsf}
\end{tabular}
\caption{Rotating cylinder: evolution for a coarse mesh of $2^4\times 2^4$ 
points and 4 adaptive refinement levels. Snapshots of the cylinder and the 
corresponding adaptive mesh:
(upper) after one time step, (lower) after 1/2 turn.}
\label{cyl4442}
\end{figure}
Figure \ref{cyl4442} represents the evolution of the rotating
cylinder on a half turn  with a coarse mesh of  $16\times 16$ points
and 4 adaptive refinement levels. We notice that the cylinder is well
represented and that the mesh points concentrate along the discontinuities.

\subsection{The two-stream instability}
\label{sec:TSI}
\begin{figure}
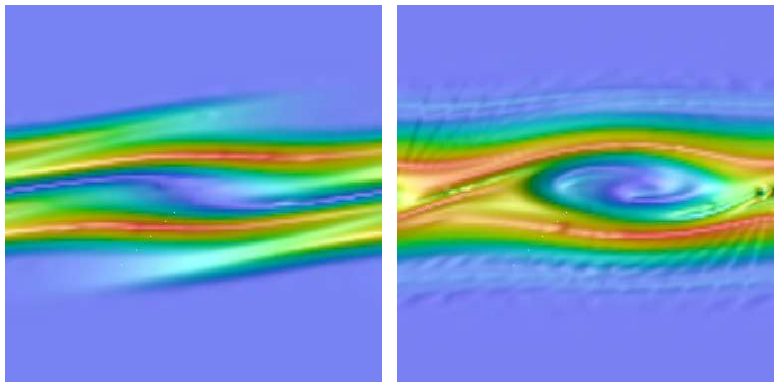

\centerline{
\begin{tabular}{ccc}
\includegraphics[width=5cm,angle=-90]{fxv.5.epsf} &\hspace{1cm}&
\includegraphics[width=5cm,angle=-90]{fxv.30.epsf}
\end{tabular}}
\caption{Two stream instability for  a coarse mesh of $2^4\times 2^4$, and 3
adaptive refinement levels,
(left) at time $t=5 \omega_p^{-1}$, (right) at time $t=30 \omega_p^{-1}$.}
\label{t5_3}
\end{figure}

We consider two streams symmetric with respect to $v=0$ and represented by
the initial distribution function
$$
f(0,x,v) = \frac{1}{\sqrt{2\pi}}v^2\,\exp(-v^2/2)(1+ \alpha \,\cos(k_0\,x)),
$$
with $\alpha=0.25$, $k_0=0.5$, and $L=2\,\pi/k_0$. We use a maximum of
$N_x=128$ points in the $x$ direction, and $N_v=128$ points in the
$v$ direction with $v_{max}=7$, and a time step $\Delta
t=1/8$. The solution varies first very slowly and then fine scales are
generated.
Between times of around $t \simeq 20 \, \,\omega_p^{-1}$ and $t
\simeq 40 \, \,\omega_p^{-1}$, the instability increases rapidly and a hole
appears in the middle of the computational domain. After 
$t = 45 \,\,\omega_p^{-1}$
until the end of the simulation, particles inside the hole are trapped.
On figure \ref{t5_3} we show a snapshot of the distribution function at times
$t = 5 \,\,\omega_p^{-1}$ and $t = 30 \,\,\omega_p^{-1}$ for a coarse mesh of 
$16\times 16$ points and 3 levels of refinement. The adaptive method
reproduces well the results obtained in the non adaptive case.
\section{Conclusion}
In this paper we have described a new method for the
numerical resolution of the Vlasov equation using an adaptive
mesh of phase-space. The adaptive algorithm is based on a 
multiresolution analysis. It performs qualitatively well.
However, there is a large overhead due to the handling of the adaptive
mesh which has not been optimized yet.
The performance of the code needs to be improved before we
can recommend this technique for actual computations. We are currently
working on optimizing the code and trying different kinds of wavelets,
as well as obtaining error estimates for the adaptive method.

\end{document}